\documentclass{amsart}
\usepackage{amsfonts,amscd}  

\newtheorem{claim}{}[section]
\newtheorem{theorem}[claim]{Theorem}
\newtheorem{lemma}[claim]{Lemma}
\newtheorem{proposition}[claim]{Proposition}
\newtheorem{corollary}[claim]{Corollary}

\newtheorem{example}[claim]{Example}

\def\proclaim #1. #2\par{\medbreak
\noindent{\bf#1.\enspace}{\sl#2}\par\medbreak}
\makeatother
 
\DeclareMathOperator{\ca}{C^*-\text{algebra}}

\DeclareMathOperator{\cas}{C^*-\text{algebras}}  
\DeclareMathOperator{\wa}{W^*-\text{algebra}}  
\DeclareMathOperator{\was}{W^*-\text{algebras}}  
\DeclareMathOperator{\csa}{C^*-\text{subalgebra}} 
\begin{document}

\title{Complete 
 isometries into $C^*$-algebras}

\date{March 5, 2002}   
 
\author{David P. Blecher and Damon M. Hay}
\address{Department of Mathematics, University of Houston, Houston,
TX 77204-3008}
\email[David P. Blecher]{dblecher@math.uh.edu}
\email[Damon Hay]{dhay@math.uh.edu}
\thanks{This research was supported in part by grants from
the National Science Foundation}

\begin{abstract} 
We give various characterizations of into (not necessarily
onto) complete isometries between $C^*$-algebras, generalizing
a classical result of Holsztynski.   Our results
are related to a natural embedding of the noncommutative 
Shilov boundary in a second dual.   
\end{abstract}  
 
\maketitle

 
\let\text=\mbox

 
\section{Introduction.}  The classical Banach-Stone theorem 
characterizes onto linear isometries between $C(K)$ spaces,
that is, between commutative unital $\cas$.   
These results have 
been extended in many directions.
We are concerned in the present paper with linear isometries
which are not necessarily surjective (i.e. onto).  
In the case of $C(K)$ spaces
these isometries were characterized by Holsztynski
\cite{Hol}:

\begin{theorem} \label{mat}  Let $K_1, K_2$ be compact Hausdorff
spaces.   A linear map $T :  C(K_1) \rightarrow C(K_2)$ 
is an isometry if and only if $T$ is
contractive, and there exists a closed subset $E$ of $K_2$,
and two continuous functions $\gamma : E \rightarrow \mathbb{T}$
and $\varphi : E  \rightarrow K_1$, with $\varphi$
surjective, such that
$$T(f)(y) \; = \; \gamma(y) f(\varphi(y))$$
for all $y \in E$.
\end{theorem}

Here $\mathbb{T}$ is the unit circle.
Informally, the result is saying
that for any isometry $T$ 
there is a certain `part $E$ of the space $B = 
C(K_2)$ acts on',
such that $T$ `compressed to this part' has a particularly nice form,
a form which amongst other things ensures that $T$ is an isometry.
The action of $T$ on the `complementary part' plays little  
role in the fact that $T$ is an isometry.
If $T(1) = 1$, then the $\gamma$ may be omitted in the 
theorem, and then $T$ compressed to $E$ is the composition
operator $f \mapsto f \circ \varphi$, which is
an isometric homomorphism.

In the present paper we state  
noncommutative versions
 of this result.   For surjective isometries between 
$\cas$ the first `noncommutative Banach-Stone theorem'
is due to R. V. Kadison \cite{Ka}.  
We will use 
ideas of Arveson, Choi and Effros, Hamana, and Kirchberg
(although our paper is fairly self-contained).
We will consider (not necessarily 
surjective) linear 
maps $T : A \rightarrow B$ between $\cas$, and establish several 
criteria which are each equivalent to $T$ being a
linear complete isometry.

To introduce the first such criterion, we will need some 
simple notation.
Suppose that  $K,H$ are Hilbert spaces, with
closed subspaces $L, M$ respectively, that $X$
is an operator space, and
that  $T : X \rightarrow B(K,H)$.   It is clear
that the following are equivalent:
\begin{itemize}
\item [(i)]
$T(X) L \subset M$ and $T(X) L^\perp \subset
M^\perp$,
\item [(ii)]  $T(X) L \subset M$ and $T(X)^* M \subset L$,
\item [(iii)]  $T(x) P_L = P_M T(x) P_L = P_M T(x)$ for all
$x \in X$.
\end{itemize}
In this case  there are unitaries $U$ and $V$
such that $U T(\cdot) V$ is of `block diagonal' form
$$\left[ \begin{array}{ccl}
R(\cdot) & 0 \\ 0 & S(\cdot) \end{array} \right]$$
where $R : X \rightarrow B(L,M)$ and $S : X \rightarrow
B(L^\perp,M^\perp)$.   For example, $R = P_M T(\cdot) P_L$.
We have $\Vert T(x) \Vert = \max \{ \Vert R(x) \Vert ,
\Vert S(x) \Vert \}$ for $x \in X$, and similarly
for matrices.  In this case we say that $R$ is a {\em reducing
compression} of $T$.
 
A similar notion exists in the more abstract setting
of maps $T : X \rightarrow M$  from an
operator space $X$ into a $\wa$ $M$.  We say that
$R$ is a reducing 
compression of $T$, if there exist projections
$p, q \in M$ such that $R = q T(\cdot) p = T(\cdot) p =
q T(\cdot)$.   Setting $S = (1-q) T(\cdot) (1-p)$,
we have $R + S = T$, and $\Vert T(x) \Vert = \max \{ \Vert R(x) \Vert ,
\Vert S(x) \Vert \}$ for $x \in X$, and similarly
for matrices.    Representing the $\wa$ $M$ on a Hilbert
space, and choosing orthonormal bases suitably,
we may write $T(\cdot)$ as a block diagonal matrix
 $R(\cdot) \oplus S(\cdot)$ as in the
previous paragraph.
Next, if $T : X \rightarrow B$ maps into a $\ca$ $B$, we 
will say that a map $R : X \rightarrow B^{**}$ is a
reducing $\wedge$-compression 
of $T$, if $R$ is a reducing compression of
$ \; \hat{} \;  \circ T$, where $\; \hat{} 
\; $ is the canonical
1-1 *-homomorphism of $B$ into the $\wa$ $B^{**}$.  
If we take the universal representation of $B$ on 
a Hilbert space $H_u$, then $B^{**} \subset B(H_u)$, and
there are two subspaces of $H_u$ with respect to which 
$T$ has a  $2 \times 2$ matrix representation $R \oplus S$
as in the last paragraph.
 
We may now state our first characterization of 
complete isometries.  Namely a linear map  
$T : A \rightarrow B$ between $\cas$
is a complete isometry if and only if
$T$ is completely contractive and 
$T$ has a reducing 
$\wedge$-compression of the 
form $u \pi(\cdot)$, 
for a *-homomorphism $\pi$ defined on $A$,
and a certain partial isometry $u$.  
If further $T$ is unital (i.e. $T(1) = 1$), then 
$u = 1$ and 
the $S$ above is completely positive.   This is a
generalization of an old result of Choi and Effros
(7.1 in \cite{CEIOS}).   

Again we see as 
in the commutative case that it is the 1-1 *-homomorphism
`part'  which
makes $T$ a complete isometry, the $S$ plays little role.
 
Unfortunately, the $\pi$ and $S$ 
mentioned above
map into $B^{**}$ and not into $B$.  Considering
the commutative results mentioned above, one would expect
the 1-1 *-homomorphism $\pi$ to map into a 
quotient of $B$ by a closed ideal $J$ (recall that in the 
commutative case, quotients of $C(K)$ by closed ideals 
correspond to closed subsets $E$ of $K$).  
We do have such a variant of 
our theorem.  For example if $A, B$ are unital $\cas$,
then we have a characterization of 
complete isometries which is a precise noncommutative analogue of 
\ref{mat},
but in terms of a left ideal $J$ in $B$,
a *-homomorphism $\pi$ from $A$ into some sort of quotient
of $B$ involving $J$,
and a partial isometry $u$ in this quotient.    
See Theorem \ref{ma} 
(vi).  However a 
noncommutative obstacle arises: unless $J$ is a two-sided
ideal (and this is not generally possible in our situation),
then this quotient is not generally 
a $\ca$, nor is it an operator algebra.    (See \cite{Ki2}
Lemma 2.3 for additional discussion of this point.) 
Nonetheless, this is not a significant obstacle, since 
the calculations we need to do in this `uncomfortable
quotient' work out quite easily.   In Section 2,
we discuss some elementary facts about these quotients.         
Also, these quotients of $B$
contain a $\ca$ for which 
the products we need 
to consider in these quotients are well defined.
Complete isometries between $\cas$ are discussed in  Section 3. 

Our results generalize immediately to nonselfadjoint 
algebras.  While the adaption is quite trivial,
there are some interesting applications of this generalization 
which we  discuss in \cite{BL}.

We remark that our methods depend strongly on the `complete 
isometry' condition; and thus have no overlap with work on
isometries of operator algebras.
 
We end this introduction with some notation and basic facts 
which will be useful.  We write $H, K, L$ for Hilbert spaces.
Perhaps confusingly, we will use the symbols $I, J, K$ for 
ideals in a $\ca$.   All ideals are taken to be uniformly closed.
A projection on a Hilbert space, or in a
$C^*$-algebra, will mean an orthogonal projection.
If $Y$ is a subspace of $X$ we write $q_Y$ or $q$ for the
natural quotient map $X \rightarrow X/Y$.  If $S$ is a
subset of a $\ca $ $B$, then we write $C^*_B(S)$ for the 
$C^*$-subalgebra of $B$ generated by $S$.   
By `c.a.i' we mean a contractive approximate identity.  

If $X, Y$ are subspaces of a Banach algebra, we write 
$XY$ for the {\em uniform closure } of the set of finite sums 
of products of the form $x y$ for $x \in X, y \in Y$.  
For example, if $J$ is a left ideal of a $C^*$-algebra $A$,
then with this convention $J^* J$ will be
a $C^*$-algebra.  This convention extends to
three sets, thus $J J^* J = J$ for a left ideal of a $C^*$-algebra,
as is well known.

For the purposes of this paper
we define a {\em triple system}
to be a (uniformly closed) subspace $X$
of a $\ca$ such that $X X^* X \subset X$.    By `subspace' we
will allow for example spaces such as $B(K,H)$, regarded
as the `1-2-corner' of the $C^*$-algebra
$B(H \oplus K,H \oplus K)$ in the
usual way.  The important
structure on a triple system is the `triple product'
$x y^* z$.  A `triple subsystem' is a
uniformly closed vector subspace of a triple system which is closed
under this triple product.

It is well known  that $X X^* X = X$
for a triple system $X$.  Also,
it is clear that $X X^*$ and $X^* X$ are $C^*$-algebras,
which we will call the left and right $C^*$-algebras of $X$ respectively,
and $X$ is a $(X X^*) - (X^* X)$-bimodule.
A linear map
$T : X \rightarrow Y$ between triple systems is a {\em triple
morphism} if $T(x y^* z) = T(x) T(y)^* T(z)$ for all
$x,y,z \in X$.   Triple systems are operator spaces, and
 triple morphisms behave very similarly to *-homomorphisms 
between $\cas$: triple morphisms are automatically
completely contractive and have closed range.
A triple morphism is
completely isometric if it is 1-1.  The kernel of a
triple morphism on $X$ is a `triple ideal' (that is, a
uniformly closed
 $(X X^*)-(X^* X)$-subbimodule).  The quotient of a 
triple system by a triple ideal is a triple system in an
obvious way.  If one factors a triple morphism by 
its kernel one obtains a 1-1 triple morphism on the
quotient triple system.   Finally, a 
triple morphism $T : X \rightarrow Y$ between triple systems 
canonically induces a *-homomorphism $\pi : X^* X 
\rightarrow Y^* Y$ between the associated
right $\cas$, via the prescription $\pi(x^* y) = T(x)^* T(y)$.   
Similarly for the left $\cas$.
These results are all rather old, may be found in \cite{Ham3}, and
are related to results of Harris and Kaup.

It is well known that
triple systems are `the same thing'
as Hilbert $C^*$-modules, although there is a slight difference
of emphasis in the two theories.  In particular,
it is sometimes convenient to view triple systems
as {\em strong Morita equivalence bimodules}, and vice versa.
This allows us to apply the known theory of Morita equivalence
\cite{Ri2,Rieffel2}.
If $A$ and $B$ are $\cas$, then an equivalence $A-B$-bimodule $Y$
is an $A-B$-bimodule which is also both a full left and a full
right Hilbert $C^*$-module (over $A$ and $B$ respectively),
such that the left and right inner products are compatible:
$$[x,y] z = x \langle  y , z \rangle$$
for all $x, y, z \in Y$.  It is easy to see (by looking
at the so-called `linking $\ca$') that such $Y$ 
is a triple system, and $Y^* Y \cong B$ and $Y Y^* \cong A$
*-isomorphically.   Conversely, any triple system is an
equivalence  $(YY^*)-(Y^*Y)$-bimodule.   Triple ideals 
in an equivalence $A-B$-bimodule are precisely the closed
$A-B$-submodules.

If $X$ is a `unital operator space' (i.e. an operator 
space containing a fixed element $e$, such that there 
exist complete isometries $T : X \rightarrow B(H) $ with 
$T(e) = I_H$), then there exists a $C^*$-envelope of $X$,
namely a pair $(B,j)$ consisting of a unital $\ca$ $B$
and a `unital' complete isometry $j : X \rightarrow 
B$ whose range generates $B$ as a $\ca$, with the following
universal property:  For any other pair $(A,i)$ consisting of a unital $\ca$
and a `unital' complete isometry $i : X \rightarrow A$
whose range generates $A$ as a $\ca$, there exists a (necessarily 
unique, unital, and surjective) *-homomorphism $\pi : 
A \rightarrow B$ such that $\pi \circ i = j$.  
This we call the Arveson-Hamana theorem \cite{Arv1,Ham},
and as customary
we write $C^*_e(X)$ for $B$ or $(B,j)$ (it is essentially
unique, by the universal property).  See also 
\cite{BShi} for a more thorough discussion of this, and of the
next concept too.  

More generally, if $X$ is an operator space, then there exists a 
{\em triple envelope} of $X$,
namely a pair $(Z,j)$ consisting of a triple system $Z$ 
and a linear complete isometry $j : X \rightarrow
Z$ whose range generates $Z$ as a triple system
(that is there exists no nontrivial closed 
triple subsystem containing $j(X)$), with the following
universal property:  For any other pair $(W,i)$ consisting of a triple system 
and a  complete isometry $i : X \rightarrow W$
whose range generates $W$ as a triple system, there exists a (necessarily
unique and necessarily surjective) triple morphism $\pi :
W \rightarrow Z$ such that $\pi \circ i = j$.
This theorem 
(dating to the late `80's) is due to Hamana \cite{Ham3},
and we write ${\mathcal T}(X)$ for $Z$ or $(Z,j)$ (it is essentially
unique, by the universal property).  
It is easy to show that 
a $C^*$-envelope of a unital operator space is a triple envelope
too.
 
A partial isometry of a triple system $Z$ will be an
$u \in Z$ with $u u^* u = u$.  We say that $u$ is an 
isometry of $Z$ if $u^* u$ is the identity of the right $\ca$
of $Z$; or equivalently, if $z u^* u = z$ for all 
$z \in Z$.  A similar definition defines a
coisometry of $Z$; and a unitary of $Z$ is a 
coisometry which is also an isometry of $Z$.  These coincide 
with the usual definitions if $Z$ is a $\ca$.

\begin{corollary} \label{esy}  If $T : Z \rightarrow W$ is a
(not necessarily
surjective) complete isometry between triple systems
or $\cas$, then if  
$T(z)$ is a partial isometry (resp. isometry,
coisometry, unitary) of $W$ then 
$z$ is a partial isometry (resp. isometry,
coisometry, unitary) of $Z$.
\end{corollary}

\begin{proof}  Let $W'$ be the triple subsystem of $W$ 
generated by $T(Z)$.  By the last mentioned theorem, 
there is a surjective triple morphism $\tau : W' \rightarrow Z$
taking $T(z)$ to $z$.   The result is now clear.
 \end{proof} 

A similar result holds if $Z$ is merely an operator space
(resp. unital operator space), but then the conclusion
`of $Z$' needs to be replaced by `of ${\mathcal T}(Z)$ (resp.
`of $C^*_e(Z)$).  For example,
such considerations allow one to say a little more 
about `unital operator spaces'.  Suppose that $(X,e)$ is a 
unital operator space in the sense above, and suppose that 
$v$ is another element in $X$ for which there exists 
a complete isometry $S : X \rightarrow B(K)$ with $S(v) 
= I_K$.    Let $(C^*_e(X),j)$ be a
$C^*$-envelope of $(X,e)$.  Then the triple subsystem $W$ of 
$B(K)$ generated by $S(X)$
is a unital $\ca$, and by Hamana's theorem mentioned above, there
exists a surjective triple morphism $\tau : W \rightarrow 
C^*_e(X)$ taking $I_K$ to $v$.   Thus as in the proof above,
$j(v)$ is a unitary of $C^*_e(X)$, and it is not hard to 
see that $(C^*_e(X),j(v)^*j(\cdot))$
is a $C^*$-envelope of $(X,v)$.

\section{Preliminary results.}    

The first several of the following Lemma's are standard in functional
analysis.

\begin{lemma} \label{basdu} 
 Let $X$ be a closed linear subspace of a Banach space $Y$.
Then 
\begin{itemize}
\item [(1)]  As subsets of $Y^{**}$ we have
$\overline{(\hat{X})}^{w*} = X^{\perp \perp}$,
\item [(2)]   The second dual of the inclusion map 
$i : X \rightarrow Y$ is a complete isometry of 
$X^{**}$ onto $X^{\perp \perp}$.   Thus 
$X^{**} \cong X^{\perp \perp}$ isometrically, via
this canonical isometry.  
\item [(3)]   $(Y/X)^{**} \cong Y^{**}/X^{\perp \perp}$
isometrically, and weak*-weak*-homeomorphically,
via a canonical map, namely the `transpose' of the canonical
isomorphism $(Y/X)^{*} \rightarrow X^{\perp}$.
This is the same as the map obtained from
$q^{**}$, where $q : Y \rightarrow Y/X$ is
the canonical quotient map, by factoring out
Ker $q^{**} = X^{\perp \perp}$.
\item [(4)]  $\hat{Y} \cap X^{\perp \perp} = \hat{X}$.
 \end{itemize}    
\end{lemma} 

\begin{proof}  The proofs of (1), (2) and (3) are in any book
on functional analysis.   

(4): We have 
$\hat{X} \subset X^{\perp \perp} \subset Y^{**}$.  
If $\hat{y} \in \hat{Y}$ but $\hat{y} \notin \hat{X}$,
choose $\phi \in Y^*$ such that $\phi(X) = 0$,
but $\phi(y) \neq 0$.  If $\hat{y} \in X^{\perp \perp}$
there exists by (1) a net
$a_i \in X$ with $\hat{a_i} \rightarrow \hat{y}$
in the weak*-topology of $Y^{**}$.  Hence
$\phi(a_i) = 0 \rightarrow \phi(y) \neq 0$!
\end{proof}

Sometimes we will simply write $X$ for the subset $\hat{X}$ of
$X^{**}$, and often we will identify
$X^{**}$ with the subspace $X^{\perp \perp}$ of $Y^{**}$.

\begin{lemma} \label{leftid}  Let $J$ be a weak* closed
left ideal in a $\wa$ $M$.   Then there exists a unique projection 
$e  \in M$ with $J = Me$.    Thus $M/J \cong M(1-e)$ 
completely isometrically.   If $J$ is a 
weak* closed two-sided ideal, then $e$ is in the center $Z(M)$.   
\end{lemma}
 
\begin{lemma} \label{leftidca}  Let $J$ be a left ideal
in a $\ca$ $A$.   Then $J^{\perp \perp}$ is
a weak* closed left ideal in the $\wa$ $A^{**}$.  Thus
there exists a projection
$e  \in A^{**}$ with $J^{\perp \perp} = A^{**} e$,
and 
$$(A/J)^{**} \cong  A^{**}/J^{\perp \perp} \cong  A^{**} (1-e)$$ 
completely isometrically.  The projection $e$ may be taken to be
any weak* limit point in $A^{**}$ of a right c.a.i. for $J$.

Thus (composing with the canonical injection of
$A/J$ into its second dual), there is
a canonical completely isometric
linear map $$d_J  \; :\;
 A/J \rightarrow  A^{**} (1-e),$$
namely the map $a + J \mapsto  \hat{a} (1-e)$.

If $J$ is a  two-sided ideal, then $e$ is in the center $Z(A^{**})$.  
\end{lemma}

We will call the projection $e$ in the last lemma the {\em 
support projection of} $J$.

Note that by the last Lemma
we may view the quotient $A/J$ of a
$\ca$ by a left ideal, as a subspace of the left ideal
$A^{**} (1-e)$ of $A^{**}$.   Thus 
$A/J$ possesses a `partial
triple product' $\{ \cdot , \cdot , \cdot \}$,
namely $\{ x , y , z \} = x y^* z$ if $x,y,z \in A/J$ and 
if $x y^* z$ (viewed as a product in $A^{**} (1-e)$) is again in
$A/J$.   

Similar considerations apply for quotients by right ideals.
Also, $A/J$ possesses a natural `involution'
$A/J \rightarrow A/J^*$.  If we write this 
involution as $^*$, then $(a + J)^* = a^* + J^*$ for 
$a \in A$.    

Next we consider 
the quotient $A/(J + J^*)$.   That the sum of a
left ideal and a right ideal in a $\ca$ is norm closed
is apparently first due to Combes and Perdrizet;
Rudin showed this more generally in a Banach algebra
(if $J$ has a bounded right approximate identity).   See 
\cite{Dixon} for references, and
for amendation of a proof due to Kirchberg
that does the  Banach algebra case too.

We need a generalization of Lemma \ref{leftidca}:

\begin{lemma} \label{leftidca2}  (\cite{Ki})  Let $J$ be a left ideal
in a $\ca$ $A$, and let $K$ be a right ideal in $A$.
Then $J + K$ is norm closed in $A$.
Let $e$ be the support projection of $J$ in 
$A^{**}$ (see  Lemma \ref{leftidca}), and let $f$ be the support projection of 
$K$.
Then  the weak* closure of $J + K$ in
$A^{**}$, equals $f A^{**} + A^{**} e$.  Also,
$$\left( A/(J + K) \right)^{**} \cong  A^{**}/(J + K)^{\perp \perp} 
\cong  (1-f) A^{**} (1-e)$$
completely isometrically.   Thus (composing with the canonical injection of 
$A/(J + K)$ into its second dual), there is
a canonical completely isometric  map $$c_{J + K} \; :\;
 A/(J + K) \rightarrow (1-f) A^{**} (1-e),$$
namely the map $a + (J + K) \mapsto (1-f) \hat{a} (1-e)$.
 \end{lemma}

\begin{proof}   That $J + K$ is norm closed is
mentioned above.
We provide the simple proof of the other assertions.  Consider
the completely contractive weak* continuous projection 
$\Phi : x \mapsto (1-f) x (1-e)$ on $A^{**}$.   The kernel of 
$\Phi$ is $f A^{**} + 
A^{**} e$.   Thus $f A^{**} +
A^{**} e$ is weak* closed, and therefore it contains the
weak* closure of $J + K$.  However, it is fairly clear 
(using Lemma \ref{leftidca} and 
\ref{basdu} (1) if necessary) that 
$J + K$ is weak* dense in $f A^{**} + A^{**} e$.  Thus
Ker $\Phi$ is the weak* closure of $J + K$.   
Since $\Phi$ is a complete quotient map, the 
first displayed equations follow.   The rest is clear.
 \end{proof}  

In the last result, if
$K = J^*$, then $e = f$.
Since $(1-e) A^{**} (1-e)$ is then a $\ca$, we may regard 
the quotient $A/(J + J^*)$ above as a `partial $\ca$', and write 
$x * y$ for the partial product on $A/(J + J^*)$.  This 
is not the product one would think of first; it is not 
$q(a_1) * q(a_2) = q(a_1 a_2)$ for  most $a_1, a_2 \in A$
(here $q = q_{J + J^*}$).   
Extending this notation, 
we will also write $x * y = z$ if $x, z \in A/J$,
$y \in A/(J + J^*)$, and 
if $x y = z$ in $(1-e) A^{**} (1-e)$.   That is, the
natural module action $A^{**} (1-e) \; \times \; (1-e) A^{**} (1-e)
\rightarrow \; A^{**} (1-e)$ makes $A/J$ a `partial right
$A/(J + J^*)$-module'.  In fact  $A/J$ has an 
$A/(J + J^*)$-valued inner product.   This is because the 
canonical bilinear map $A/J^*  \; \times \; A/J 
\rightarrow A/(J + J^*)$ is easily seen to be well
defined.  Thus we may think of $A/J$ as a 
`partial right Hilbert $C^*$-module' over $A/(J + J^*)$.   
Similarly, $A/J^*$ is a `partial left 
Hilbert $C^*$-module' over $A/(J + J^*)$. 

As  noted in \cite{Ki,Ki2}, if $A$ is a unital $\ca$ then the
proof above of
\ref{leftidca2} shows that $A/(J + J^*)$  is an operator system.
Indeed it is  what Kirchberg calls a {\em $C^*$-system}, namely an
operator system whose second dual is a $\ca$.  We shall not
need this notation, but we refer the reader to his
very deep papers (see e.g. \cite{KW} and references therein)
for an astonishing converse,
and a  host of related notions, and deep consequences.
If $A$ is unital, we
 recall from \cite{Ki,Ki2} the `two-sided multiplier algebra' 
$M(A/(J + J^*))$.   If $A$ is unital, then $M(A/(J + J^*))$
is a unital subspace of $A/(J + J^*)$ which is actually a 
$\ca$, indeed which may be viewed as 
a *-subalgebra of $(1-e) A^{**} (1-e)$.   More
generally, if $A$ is
nonunital, we define $M_0(A/(J + J^*))$ to be the 
set $\{ b \in A/(J + J^*) : b * (A/(J + J^*))
\subset A/(J + J^*) \; , \; (A/(J + J^*)) * b 
\subset A/(J + J^*) \}$.   Here $*$ is the 
`partial product' on $A/(J + J^*)$ discussed 
just after the proof of \ref{leftidca2}.  Then 
$M_0(A/(J + J^*))$ is a *-subspace of 
$A/(J + J^*)$ which is also a $\ca$, indeed which may be 
viewed as a *-subalgebra of $(1-e) A^{**} (1-e)$. 
If $A$ is unital then $M_0(A/(J + J^*)) = M(A/(J + J^*))$.         
It is interesting that if  $A$ is unital then $w = a + (J + J^*)$ is in 
$M(A/(J + J^*))$ if and only if 
$w, w * w^*$ and $w^* * w$ are in $A/(J + J^*)$ (see
\cite{Ki}).   We will not use this fact though.

We say that a linear map $T : C \rightarrow A/(J + J^*)$ 
defined on an algebra $C$ is
a {\em homomorphism} if $T$ maps into
$M_0(A/(J + J^*))$, and $T$ is a homomorphism into 
this algebra.   A similar notation applies to 
*-homomorphisms into $A/(J + J^*)$. 
Finally, we say that 
$u \in A/J$ is a  {\em partial isometry ' 
of} $A/J$ if  $u^* u$ is a projection in 
$M_0(A/(J + J^*))$.   

We remark that if $A$ is commutative, then $J = J^* = 
J + J^*$, $A/(J + J^*)$ is an algebra,
and `homomorphisms' in the sense above
 into $A/(J + J^*)$ are exactly the 
algebra homomorphisms.   

There seems to be another $C^*$-subalgebra of $M_0(A/(J + J^*))$ 
which is of interest.   Namely, we let 
$M_{00}(A/(J + J^*))$ be the 
set $\{ b \in A/(J + J^*) : b * (A/J^*)
 \subset A/J^* \; , \; (A/J) * b
\subset A/J \}$.   Here $*$ are the
`partial module products' discussed
just after the proof of \ref{leftidca2}.  Then
$M_{00}(A/(J + J^*))$ is a *-subspace of
$A/(J + J^*)$ which is also a $\ca$, indeed which may be
viewed as a *-subalgebra of $(1-e) A^{**} (1-e)$.   
  
If $J$ is a left ideal and $K$ a right ideal of a
$\ca$ $A$, then we say that a linear map $T : Z \rightarrow A/(J
+ K)$
defined on a triple system $Z$ is a  {\em partial triple 
morphism} if $T$ composed with the
canonical map $c_{J + K} : 
A/(J + K) \rightarrow (1-f) A^{**}(1-e)$
discussed in \ref{leftidca2}, is a triple morphism.

The following result then
follows immediately from well known facts about
*-homomorphisms and triple morphisms mentioned in our 
introduction.

\begin{proposition} \label{ptri} Let $J$ be a left ideal in a $\ca$ $A$
let $K$ be a right ideal in $A$,
and let $T : C \rightarrow A/(J + J^*)$ (resp. $Z \rightarrow A/(J + K)$) 
be a *-homomorphism 
(resp. partial triple morphism)
defined on a $\ca$ $C$ (resp. triple system $Z$).     
Then $T$ is completely isometric if and only if $T$ is 1-1.
\end{proposition}   

\begin{lemma} \label{leftdid}  Let $B$ be a $\ca$, let 
$A$ be a $\csa$ of $B$ and let $I$ be a closed
two-sided ideal in $A$.  Set $J$ to be the closed left ideal of $B$
generated by $I$.   That is, $J = BI$. Then 
$J \cap A = I$.   The canonical map
$r : A/I \rightarrow B/J$ is a linear 1-1 complete isometry.
The canonical map $t : A/I \rightarrow B/(J + J^*)$ has 
range within $M_{00}(B/(J + J^*))$, and 
is a 1-1 *-homomorphism into this space.
\end{lemma}   

\begin{proof}  Certainly $I \subset J \cap A$, but if $z \in J \cap A$
and if $\{ e_i \}$ is a c.a.i. for $I$, then $z  e_i \rightarrow z$
(since $z \in J = BI
$).  However since  
$z \in A$ then $z  e_i \in I$.  Therefore $z \in I$.

Now it is easy to see that the canonical map $r : A/I 
\rightarrow B/J$ is a linear 1-1 complete contraction.
Note that $(A/I)^{**} \cong A^{**}/I^{**} \cong A^{**} (1_{A^{**}} -p)$ 
for a central projection $p \in A^{**}$ by \ref{leftid}.  
We claim that $J^{**} = B^{**} p$, so that
$B^{**}/J^{**} \cong B^{**} (1_{B^{**}} -p)$.   Certainly  since $p \in I^{**}
\subset J^{**}$, we have $B^{**} p \subset J^{**}$.  Conversely,
if $z \in J^{**}$ and if $z_i \in J$ converge weak* to $z$, then
$z_i = z_i p \rightarrow z p$, so that $z = zp \in B^{**} p$.

Note too that $\hat{a} p = p \hat{a}$ for all
$a \in A$, since $p$ is central in
$A^{**}$.    

Consider the sequence of canonical 
completely contractive maps $$A/I \rightarrow B/J \rightarrow 
B/(J + J^*) \overset{c_{J + J^*}}{\rightarrow}
 (1_{B^{**}} -p) B^{**} (1_{B^{**}} -p) \subset B^{**}.  $$ 
Explicitly, the composition $\theta$  of all of these maps takes 
$$(a + I) \mapsto
(a + J) \mapsto (\hat{a} + (J + J^*)) \mapsto 
(1_{B^{**}} -p) \hat{a} (1_{B^{**}} -p)
= \hat{a} (1_{B^{**}} -p).$$ 
If $\theta(a + I) = 0$ then $\hat{a} = \hat{a} p$.
This means that $\hat{a} \in \hat{A} \cap I^{**}$,
so that $a \in I$ by Lemma \ref{basdu} (4).
Thus $\theta$ is 1-1.

Since $p$ is central in
$A^{**}$, we have that  $\hat{B} (1-p) \hat{a} (1-p)
\in \hat{B} (1-p)$ and $(1-p) \hat{a} (1-p) \hat{B} \in (1-p) \hat{B}$.
Thus if we use the fact in
\ref{leftidca2}  that $c_{J + J^*}$ there is 1-1, we see
that $a + (J+J^*) \in M_0(B/(J + J^*))$.   Similarly,
$a + (J+J^*) \in M_{00}(B/(J + J^*))$.

To complete the proof of the 
theorem, it suffices to show that 
$\theta$ is a *-homomorphism.
However this is easy, for example note that 
$$(1_{B^{**}} -p) \hat{x} (1_{B^{**}} -p) \hat{y} (1_{B^{**}} -p) = 
(1_{B^{**}} -p) \widehat{xy} (1_{B^{**}} -p)$$
for $x,y \in A$.   
\end{proof} 

As one application of the last result, we may see very 
easily that the $C^*$-envelope $(C^*_e(X),j)$ of
a unital subspace or unital subalgebra of a unital 
$\ca$ $B$,
may be identified with a certain $C^*$-subalgebra
of the second dual of $B$:  
 
\begin{corollary} \label{ceisd}  Suppose that 
$X$ is a unital operator space, and that $T :
X \rightarrow B$ is a unital complete isometry into a
 unital $\ca$ $B$.  Then there is a 
projection $p \in B^{**}$, such that 
if we define $R : X \rightarrow (1-p) B^{**} (1-p)$ to be the 
map $R(x) = (1-p) \widehat{T(x)} (1-p)$ for $x \in X$,
then also $R(x) = \widehat{T(x)} (1-p) = (1-p) \widehat{T(x)} $
for $x \in X$, and 
$(C^*_{B^{**}}(R(X)),R)$ is a $C^*$-envelope of $X$.
\end{corollary}

\begin{proof}   Let $A = C^*_{B}(T(X))$.  By the universal
property of the $C^*$-envelope $(C^*_e(X),j)$ 
discussed in the introduction, there exists a surjective
*-homomorphism $\theta : A \rightarrow C^*_e(X)$
such that $\theta(T(x)) = j(x)$ for all $x \in X$.  Thus
if $I = $ Ker$\; \theta$, then $A/I$ is a $\ca$,
and there is a *-isomorphism $\psi :
C^*_e(X) \rightarrow A/I$ with $\psi(j(x)) = T(x) + I$
for all $x \in X$.    
If $p$ is the  support projection of
$I$ in the center of $A^{**}$, then we may view $p$ as
a projection in $B^{**}$.    
Also, there is a 1-1 *-homomorphism
$A/I \rightarrow (1-p) B^{**} (1-p)$
by the proof of Lemma \ref{leftdid}.   This map 
takes $c + I$ to $\hat{c} (1-p)$.   Composing 
this map with $\psi$ gives a 1-1 *-homomorphism
$C^*_e(X) \rightarrow (1-p) B^{**} (1-p)$ which takes
$j(x)$ to $(1-p) \widehat{T(x)}  (1-p)$. 
\end{proof}

\begin{lemma} \label{leftdidtro}  Let $B$ be a $\ca$, let
$Z$ be a uniformly closed triple subsystem of $B$ and let $N$ be a 
triple ideal in $Z$.  Set $J$ (resp. $K$)
to be the closed left (resp. right) ideal of $B$
generated by $N$.   That is, $J = BN$ and $K = NB$. Then 
$J \cap Z = K \cap Z = N$.   The canonical map 
$r : Z/N \rightarrow B/J$ (resp. $s : Z/N \rightarrow B/(J + K)$)
taking 
$z + N$ to $z + J$ (resp. 
$z + (J + K)$) are linear complete isometries. 
Indeed $r$ and $s$ are completely isometric partial triple morphisms.
\end{lemma}
 
\begin{proof}  Certainly $N \subset J \cap Z$, but if $z \in J \cap Z$
and if $\{ e_i \}$ is a c.a.i. for $I = N^* N$, then $z  e_i \rightarrow z$
(since $z \in J = BN$).  However since 
$z \in Z$ then $z  e_i \in N$.  Therefore $z \in N$.
Similarly, $K \cap Z = N$.
Also the canonical map $r : Z/N 
\rightarrow B/J$ is a linear 1-1 complete contraction.  
Set $I = N^* N$ and $D = Z^* Z$.   We claim that 
$J = BI$.   Since $B N^* N \subset B N$, we have                 
$BI \subset J$.   On the other hand, since $N$ is
a 
triple system, we have $N = N N^* N$.
Thus  $BN = B N N^* N \subset B N^* N$, so that $J \subset BI$.  
   This proves the claim.
We are also put in the situation of the previous Lemma, so that
there is a unique central projection $p$ in $D^{**}$ with 
$D^{**} p = I^{**}$ and $B^{**} p = J^{**}$.

The map we are interested in from $Z/N \rightarrow B^{**} (1-p)$ is 
precisely the map $(z + N) \mapsto \hat{z} (1-p)$.   Since 
$p$ commutes with $D^{**}$ we have for $z_1,z_2,z_3 \in Z$ that
$$\hat{z_1} \hat{z_2}^* \hat{z_3} (1-p) = 
\hat{z_1} (1-p) (1-p) \hat{z_2}^* \hat{z_3} (1-p)$$
which says that $d_J \circ r$ is a triple morphism.  Thus
$r$ is a partial triple morphism.    Since $r$
is 1-1 it is completely isometric by Proposition
\ref{ptri}.  To see that $s$ is 1-1, or equivalently, 
completely isometric we will need some information from the 
following Lemma.  Namely, if $q$ is as in the 
following Lemma, then it is easy to see (by symmetry)
that $q$ is the support projection for $K = N B$;
and that the completely isometric map 
$(z + N) \mapsto \hat{z} (1-p)$ discussed above from
$Z/N \rightarrow B^{**} (1-p)$ actually maps into
$(1-q)  B^{**} (1-p)$, and that this map factors 
completely contractively through $s$.  That is, this map
is the composition of $s$ and $c_{J + K}$ mentioned in 
\ref{leftidca2}.  Thus $s$ is a complete isometry. 
\end{proof}

Note that the ideal $K$ above is `retrievable' from 
$J$ and $B$ and $Z$.  Namely, $K = (J \cap Z) B$.
 
\begin{lemma}  \label{moce}  Let $B, Z, N$ be as in the previous 
lemma.  Associated with the support projection $p \in B^{**}$ for
$N^* N$ (mentioned 
in the previous proof), there is a projection $q \in B^{**}$,
such that $q z = q z p = z p$ for all $z \in Z$.   Indeed 
such $q$ may be taken to 
be the support projection for $N N^*$, viewed as 
a projection in $B^{**}$.
\end{lemma}
 
\begin{proof}   This will be clear to those familiar  with
the Morita equivalence of $\was$ \cite{Rieffel2}, but for the readers 
convenience we sketch a more or less elementary
and complete proof.   In the proof 
above we have $D = Z^* Z$; set $C$ to be the `left $\ca$'
$Z Z^*$ of $Z$.  Then 
it is well known from $C^*$-module theory that  
we may pick a c.a.i. $\{ e_\alpha \}$ in $C$ with terms 
$e_\alpha$ each of the 
form $\sum_{k=1}^n z_k z_k^*$.  Pick $z \in Z$.
It is an elementary 
item from $C^*$-module theory that $e_\alpha z \rightarrow z$.
Consider the associated 
net with terms $\sum_{k=1}^n z_k p z_k^*$.  
This is a net in $Ball(B^{**})$ (since $z_k p z_k^* \leq
z_k z_k^*$ for each $k$).
Let $q$ be a weak* limit point of this net
in $Ball(B^{**})$.  
We claim that $q \in M$, where $M$ is the $\wa$ which
is the weak* closure in $B^{**}$ of $N N^*$.
To see this note that $p$ was chosen in $I^{**}$, and therefore 
lies in the  weak* closure in $B^{**}$ of $N^* N$.
Thus $z_k p z_k^*$ lies in the  weak* closure in $B^{**}$ of $
Z N^* N Z^*$.
Since $N$ is a left $Z Z^*$-module, and since $N = N N^* N$
as noted in the introduction, we have $Z N^* = Z N^* N N^*
\subset N N^*$. Thus $Z N^* N Z^* \subset N N^*$.
This proves the claim.      

Since 
$$\sum_{k=1}^n z_k p z_k^* z = \sum_{k=1}^n z_k z_k^* z p 
= \sum_{k=1}^n z_k p z_k^* z p$$
in the limit we obtain $q z = z p = q z p$.  Thus
$q^2 z = q z p = q z$.   Hence $q^2 x = q x$ for all
$x \in C$, and by a weak* approximation argument, for all
$x \in M$.  Taking $x = 1_M$ shows that $q^2 = q$,
so that $q$ is a 
(contractive and therefore
selfadjoint) projection in $M$.  

Finally note that if $z, w \in Z$, then $q z w^* = z p w^*$.
But we saw above that $z p w^*$ lies in the weak* closure in $B^{**}$ of
$N N^*$.  That is, $q z w^* \in M$.  By a weak* approximation argument, 
$q R \subset M$, where $R$ is the $\wa$ which
is the weak* closure in $B^{**}$ of $Z Z^*$.  Conversely,
if $z_1, z_2, z_3, z_4 \in N$, then $q z_1 z_2^* z_3 z_4^*
= z_1 p z_2^* z_3 z_4^* = z_1 z_2^* z_3 z_4^*$,
since $p$ 
is the identity of $I^{**}$.  Since $
N N^* N N^* = N N^*$, we have $q x = x$ for all 
$x \in N N^*$, and therefore also
for all $x \in M$.  Thus $q R = M$.  
\end{proof}
 
Variants of the following lemma are no doubt well known:

\begin{lemma}  \label{maw}  Let $T : A \rightarrow Z$ be a
linear map from a $\ca$ onto a triple subsystem $Z$ of
a $\wa$ $M$.   Then $T$ is a triple morphism if and only if 
there exists a surjective *-homomorphism 
$\theta : A \rightarrow Z^* Z$, and a partial isometry 
$u \in M$, such that the projection $u^* u$ is the identity
of the $W^*$-algebra $\overline{\text{Ran} \; \theta}^{weak*}$, 
and such that $$  T(a) = u  \theta(a)
\; \; \; \; \; \; \; \; \& \; \; \;  \; \; \; \; \; \theta(a)  = u^* T(a)$$
for all $a \in A$.
 Moreover $\theta(a_1^* a_2) = T(a_1)^* T(a_2)$ for all
$a_1, a_2 \in A$; and if $T$ is 1-1, then so is 
$\theta$.  If $A$ is unital, then 
we may choose $u$ in the equivalence above to be
$T(1)$, and this is a unitary of
$Z$, and also $Z = u^* (Z^* Z)$.  If $A$ is nonunital then 
we may choose  $u$ to be any 
weak* limit point of $\{  T(e_\alpha) \}$, for any c.a.i. 
$\{ e_\alpha \}$ of $A$. \end{lemma}

\begin{proof}  (Sketch.)  If we define 
 $\theta(a_1^* a_2) = T(a_1)^* T(a_2)$ then 
as stated in the introduction,
$\theta$ is a well 
defined *-homomorphism $A \rightarrow Z^* Z$.
We have that $\theta$ is certainly onto since it 
may be clearly seen to have dense range.  If $A$ is unital 
then define $u = T(1)$, and the result is 
evident.    If $A$ is nonunital, then {\em define}
 $u$ to be any
weak* limit point of $\{  T(e_\alpha) \}$, for any c.a.i.
$\{ e_\alpha \}$ of $A$.   Then 
$$u \theta(a_1^* a_2) = 
\lim  T(e_\alpha) T(a_1)^* T(a_2) = 
\lim T(e_\alpha a_1^* a_2) = T(a_1^* a_2)$$
for $a_1, a_2 \in A$.   Thus $T = u \theta(\cdot)$.
Also, $u^* T(a) = \lim  T(e_\alpha)^* T(a)
= \lim \theta(e_\alpha^* a) = \theta(a)$.
Therefore $u^* u b = b$ for all $b \in $ Ran $\theta$,
and consequently, for all $b$ in the weak* closure
$N$ of Ran $\theta$.   It is easy to see that $u^* u$ 
is in $N$ too, so that $u^* u$ is the identity of
$N$.    
\end{proof}  

Lemmas \ref{leftid},  \ref{leftidca}, 
\ref{leftidca2}, \ref{leftdidtro} and \ref{moce} 
generalize without much change 
to `ideals' in 
triple systems,
or equivalently, to submodules of `equivalence bimodules'.
It is folklore known to rather few 
that if $Y$ is an equivalence $A-B$-bimodule
then $Y^{**}$ is an $W^*$-equivalence 
$A^{**}-B^{**}$-bimodule (see e.g.
5.8 in \cite{BShi}, \cite{Mg},\cite{EOR}).   
Weak* closed right $B^{**}$-submodules
of $Y^{**}$ are exactly the submodules of the form 
$f Y^{**}$, for a projection $f$ in $A^{**}$.
Norm closed right $B$-submodules of $Y$, are precisely 
the subspaces $R$ of $Y$ such that $R^{\perp \perp} = f
 Y^{**}$, for a projection $f$
in  $A^{**}$.   These facts may be found for example proved in
6.6 of \cite{BEZ}.  We will say that 
$f$ is the support projection of $R$.   As in Lemma
\ref{leftidca} we have $(Y/R)^{**} \cong 
Y^{**}/R^{\perp \perp} \cong (1-f) Y^{**}$ completely isometrically.
Similar results hold for left $A$-submodules $L$ of $Y$; there 
exists a `support projection' $e \in B^{**}$.  If $N$ 
is an $A-B$-submodule of $Y$, then by following the 
proof of \ref{leftdidtro} and \ref{moce}, one sees that
the left support projection
$f$ of $N$  is in
the center of $A^{**}$, the right support projection $e$ is in
the center of $B^{**}$, and 
$f \hat{y} = \hat{y} e = f \hat{y} e$ for all $y \in Y$.  

There is also a
generalization of \ref{leftidca2} to left and right 
submodules of an equivalence $A-B$-bimodule,
which is
more or less identical 
to \ref{leftidca2}.   Since we shall not explicitly 
need this it is therefore omitted, although it 
is useful in placing the proof of the
next theorem in context.   
 
On asking Kirchberg, we 
ascertained that 
the results such as the next one are known to him,
but at this point it does not appear in the literature
as far as we know.  This next result
is the `nonunital' version of \ref{ceisd}:

\begin{theorem}  \label{treisd}  Let 
$T : X \rightarrow Y$ be a complete isometry from an 
operator space into 
an equivalence $A-B$-bimodule $Y$.
Then there are projections $e \in B^{**}$ and $f \in A^{**}$,
such that $$\widehat{T(x)} e = f \widehat{T(x)} = f 
\widehat{T(x)} e$$ for all
$x \in X$, and such that if we define 
$i : X \rightarrow
Y^{**}$ by $i(x) = (1-f) \widehat{T(x)} = \widehat{T(x)} (1-e) = 
(1-f) \widehat{T(x)} (1-e)$,
and set $W$ to be the triple subsystem of $Y^{**}$ generated by
$i(X)$, then $(W,i)$ is a triple envelope of $X$.
\end{theorem}

Thus if $X$ is a subspace of a triple system $Y$, then 
we may identify the triple envelope of $X$ with  a 
certain triple subsystem of $Y^{**}$.

\begin{proof}   Let $Z$ be the  triple subsystem of $Y$
generated by 
$T(X)$, and let $({\mathcal T}(X),j)$ be a 
triple envelope of $X$.  By the universal 
property of ${\mathcal T}(X)$
discussed in our introduction, there exists a surjective
triple morphism $\mu : Z \rightarrow {\mathcal T}(X)$
such that 
$\mu(T(x)) = j(x)$ for all $x \in X$.  Thus
if $N = $ Ker $ \; \mu$, then $Z/N$ is a triple system,
and there is a triple isomorphism $\xi : 
{\mathcal T}(X) \rightarrow Z/N$ with $\xi(j(x)) = 
T(x) + N$
for all $x \in X$.    
Let  $e$ and $f$ be the right and left support projections of
the sub-bimodule $N$ of the triple system 
$Z$.   So $e \in (Z^* Z)^{**}$, for example,
and we have exactly
as in \ref{moce}, that $\hat{z} e = f
\hat{z}  = f \hat{z}  e$ for all $z \in Z$.
We will view $e$ and $f$ as 
projections in $A^{**}$ and $B^{**}$
respectively.   
Consider the map $\rho :
Z/N \rightarrow (1-f) Y^{**} (1-e)$, defined 
by $\rho(z + N) = \hat{z} (1-e)$.
If $\rho(z + N) = 0$, then $\hat{z} = \hat{z} e$.
Thus $$\hat{z} \in \hat{Z} \cap Z^{**} e 
= \hat{Z} \cap N^{\perp \perp} = \hat{N}$$
using \ref{basdu} (4).  Therefore $z \in N$, so that
$\rho$ is 1-1.   
Since $\hat{z} e = f
\hat{z}  = f \hat{z}  e$ for all $z \in Z$, 
we have that $$\widehat{x y^* z} (1-e) = 
\hat{x} (1-e) \hat{y}^* \hat{z} (1-e)$$ for $x,y,z \in Z$.   
Thus $\rho$ is a 1-1 triple morphism.
The composition $\rho \circ \xi$ is therefore a
1-1 triple morphism of ${\mathcal T}(X)$ onto $W$,
which maps $j(x)$ to $i(x)$, for each $x \in X$.
This completes the proof.  \end{proof}

The following  is a `nonsurjective Banach-Stone theorem' 
for triple systems:

\begin{corollary} \label{bsttr}  If $T : 
Z \rightarrow Y$ is a (not necessarily 
surjective) complete isometry between 
triple systems, then there exist projections
$e$ and $f$ in the second duals of 
the $\cas$ $Y^*Y$ and $YY^*$ respectively, such that 
$$(1-f) \widehat{T(z)} = \widehat{T(z)} (1-e) = (1-f) 
\widehat{T(z)} (1-e)$$
for all $z \in Z$, and that these expressions
define a 1-1 triple morphism from $Z$ into $Y^{**}$.  
\end{corollary} 
 
\begin{proof}   This follows from the previous 
result, together with the fact that if $(W,i) $ is
a triple envelope of a triple system, then necessarily 
$i$ is a triple isomorphism.   This last fact is clear from the 
universal property of the triple envelope.
\end{proof}

A similar corollary, but for unital complete isometries
between $\cas$, follows from \ref{ceisd}.  We will say more on
this case in \ref{uca}.

The last result says that if
$T  : Z \rightarrow Y$ is a  complete isometry between
triple systems, then there exists a `reducing $\wedge$-compression of $T$'
which  is a triple morphism.   Indeed with a little thought
one may see that if one takes the `universal representation' 
of the triple system $Y \subset B(K_u,H_u)$ (derived in an
obvious way from the universal representation of the `linking 
$\ca$' of $Y$), then there are projections $e$ and $f$ on 
$K_u$ and $H_u$ respectively, such that $(1-f) T(\cdot) (1-e)$ 
is a reducing compression of $T$ in the sense of the introduction,
and is a 1-1 triple morphism.  
Since we will not explicitly use this we omit the details.

\section{Banach-Stone theorems
for $\cas$.}

In the following we will often identify $B$ with the subset $\hat{B}$ of
$B^{**}$.  In particular if $b \in B$ and $p$ is a
projection in $B^{**}$, then  we will often write $b (1-p)$ for 
the expression $\hat{b} (1-p)$.
 
\begin{theorem}  \label{ma}  
Let $T : A \rightarrow B$ be a completely contractive linear map
between $C^*$-algebras.  Then the following are equivalent:
\begin{itemize}
\item [(i)]   $T$ is a complete isometry,
\item [(ii)]  there exists a $C^*$-subalgebra $D$ of $B$,  
a closed two-sided ideal $I$ in $D$, a
*-isomorphism $\pi : A \rightarrow D/I$, and a
partial isometry $u \in B^{**}$,
such that 
$$\widehat{T(a)} (1-p) = u \pi(a)$$
 and such that $u^* u \pi(a) 
= \pi(a)$, for all $a \in A$,
where $p$ is the support projection for $I$
in $D^{**}$, viewed as an
element of $B^{**}$.  We are viewing $\pi(A)$
in the equations above as a subset of 
$(1-p) B^{**} (1-p)$ via the identifications
$D/I 
\subset D^{**}/I^{**} =  (1-p) D^{**} (1-p) 
\subset (1-p) B^{**} (1-p)$,   
  \item [(iii)]  
there is a left ideal $J$ and a right ideal $K$ of $B$, such that
$q_{J + K} \circ T$ is a 1-1 partial triple morphism
from  $A$ to $B/(J + K)$,
\item [(iv)]  
$T$ possesses a reducing $\wedge$-compression 
which is a 1-1 triple morphism,
\item [(v)]  there is a projection
$p$ in $B^{**}$, a partial isometry $u$ in $B^{**} (1-p)$,
and a 1-1 *-homomorphism $\theta : A \rightarrow (1-p) B^{**} (1-p)$,
such that $$\widehat{T(\cdot)} (1-p) =
u \theta(\cdot)$$
and such that $u^* u \theta(\cdot) = \theta(\cdot)$.
\end{itemize}
Further, if $A$ is unital, then these are equivalent to
\begin{itemize}   \item [(vi)]  
there is a left ideal $J$ of $B$, a 
1-1 *-homomorphism $\pi : A \rightarrow M_{00}(B/(J + J^*))$,
and a partial isometry $u$ in $B/J$
such that $$q_J(T(a)) = u \pi(a)
\; \; \; \; \; \; \; \; \& \; \; \;
 \; \; \; \; \; \pi(a)  = u^* q_J(T(a))$$
for all $a \in A$, where 
 $q_J$ is the canonical quotient map
$B \rightarrow B/J$
(Notations as in the first and fifth paragraphs after
\ref{leftidca2}).
\item [(vii)]  there exists a $C^*$-subalgebra $D$ of $B$
containing $T(A)^* T(A)$,
a closed two-sided ideal $I$ in $D$, and a
*-isomorphism $\pi : A \rightarrow D/I$,
such that
$$q_N(T(a)) \;  = \;  U \; \; \pi(a)$$
for all $a \in A$, where $q$ is the canonical
 quotient triple morphism from the triple subsystem
$Z = T(A) D$ of $B$, onto the quotient of $Z$ by the
triple ideal $N = Z I$, and where $U$ is a 
unitary of the triple system $Z/N$.   
\end{itemize}
The $u$ in (vi) (resp. $U$ in (vii)) may be taken to 
be $q_J(T(1))$ (resp. $q_N(T(1))$).   
The projection $p \in B^{**}$ above in (ii) and 
(iv) may be chosen to commute
with (the image in $B^{**}$ of) $T(A)^* T(A)$.   One may take 
either $J$ or $K$ to be $\{ 0 \}$ in (iii) if one wishes.
 
If $A$ is nonunital then the $u$ in (ii) and (v)
may be taken to be a 
weak* limit point in $B^{**}$ of $\{ \widehat{T(e_\alpha)} (1-p) \}$,
for some 
c.a.i.  $\{ e_\alpha \}$ of $A$.
\end{theorem}

\begin{proof}   We leave it as an exercise for the reader
that any one of (ii)-(vi) imply (i).
For example, given (v),
then we have for $a \in A$ that
$$\Vert T(a) \Vert
\geq  \Vert T(a) (1-p) \Vert = \Vert u \theta(a) \Vert
\geq \Vert u^* u \theta(a) \Vert 
= \Vert \theta(a) \Vert = \Vert a \Vert.$$
Thus $T$ is an isometry.  A similar calculation shows that $T$ is a 
complete isometry, giving (i).
 
Henceforth in the proof we will assume that 
$T$ is a complete isometry.   Then (iv) follows 
immediately from \ref{bsttr} (or one may see it from the 
proof below).
To see (iii), we follow the proof of 
\ref{bsttr} 
and \ref{treisd}, 
letting $Z$ be the triple system in $B$ generated by 
$T(A)$, and setting $N$ to be the  kernel
of the 
surjective triple morphism
$\rho : Z \rightarrow A$ such that 
$\rho(T(a)) = a$ for all $a \in A$.   
Then $A$ is triple isomorphic to $Z/N$ via 
the map $\xi : a \rightarrow \; T(a) + N$.   
Apply \ref{leftdidtro} to $Z/N$.
We obtain a 1-1 partial triple morphism
$\psi : Z/N \rightarrow B/(J + K)$, and by composition, a 
1-1 partial triple morphism $\theta : A \rightarrow  B/(J + K)$.
We have $\theta(a) = T(a) + (J + K)$.   This implies (iii).

To see (ii), consider again the triple isomorphism 
$\xi : A \rightarrow Z/N$ above.  Via the proof 
of Lemma \ref{leftdidtro}
we may view $Z/N$ as a triple subsystem $W
= \hat{Z}(1-p)$ of $B^{**}(1-p)$,
via the map $z + N \mapsto \hat{z} (1-p)$.  We saw in the proof 
of \ref{leftdidtro} that $p$ commutes with the image of 
$Z^* Z$ in $B^{**}$.  Thus 
$$T(a_1)^* T(a_2) (1-p) = (1-p) T(a_1)^* T(a_2) = (1-p) T(a_1)^* T(a_2)
(1-p)$$
for all $a_1, a_2 \in A$.    Let $D = Z^* Z$, $I = N^* N$.   
By Lemma \ref{maw} we have $T(a)(1-p) = u \pi(a)$, for a 
partial isometry $u \in B^{**}(1-p)$,
and a *-isomorphism $\pi'$ from $A$ onto
$W^* W$, which may be defined to be   
$\pi'(a) = u^* T(a)(1-p)$; or alternatively,
$\pi'(a_1^* a_2) = (1-p) T(a_1)^* T(a_2) (1-p)$
if $a_1, a_2 \in A$.  We saw in 
Lemma \ref{maw} that $u$ may be taken to be a weak* limit 
point of $\{ T(e_\alpha)(1-p) \}$ in $B^{**}$, for any
c.a.i.  $\{ e_\alpha \}$ of $A$.     
From the basic theory 
of Morita equivalence of quotients
(see \cite{Rif2}), or the basic theory
of quotients of triple systems \cite{Ham3},
 $W^* W \cong (Z^* Z)/(N^* N) = D/I$
*-isomorphically, via the natural map.  
Write $\sigma$ for the *-isomorphism $A \rightarrow D/I$
corresponding to $\pi'$;
it is clear that $\sigma(a_1^* a_2) = 
T(a_1)^* T(a_2) + I$, for $a_1, a_2 \in A$.
This gives (ii), and also (v).

Tracing the proof above of (v), we see that
$\theta$ maps into
$(1-p) Z^* Z (1-p) \subset 
c_{J + J^*}(B/(J + J^*))$.
Set $\pi = c_{J + J^*}^{-1} \circ \pi' = t \circ \sigma$
(where $t$ is as in Lemma \ref{leftdid}), and $u = q_J(T(1))$,
to obtain (vi).         

If $Z, N, D, I, \xi$ are as 
in our proof of (ii), and if one applies Lemma
\ref{maw} to $\xi$, and uses the fact (also essentially
used in the proof of (ii)) that the right
$C^*$-algebra of the triple system $Z/N$ is *-isomorphic
to $D/I$, one obtains (vii).
\end{proof}

From (iv) of the theorem above we separate the following
fact:
$$T(a) T(b)^* T(c) (1-p) = T(a b^* c) (1-p)$$
for all $a, b, c \in A$.

\vspace{3 mm}

{\bf Remarks.}   1).   In the proof of (vii) we have $D = Z^* Z$ and 
$I = N^* N$.   By the theory of quotient
Morita equivalences \cite{Rif2} or quotient triple
systems, $Z/N$ is a full right $C^*$-module over
$D/I$, so that the product in the displayed equation in
(vii) makes perfect sense.

2.  It is not mentioned explicitly in the theorem, but looking at
the details of the proof, one sees that there
is a canonical projection $q \in B^{**}$ such that $(1-q) \widehat{T(a)}
= \widehat{T(a)} (1-p)$ for all $a$.

3.  It is easy to recover Holstynski's theorem \ref{mat}
directly from (vi) of the last theorem.   We leave this as an 
exercise for the reader (the point is that firstly,
because the isometry maps between subspaces of
$C(K)$ spaces, it is completely
isometric (\cite{Pau} Theorem 3.8); and secondly in the 
commutative case left ideals are two-sided, so that 
$J = J^* = J + J^*$.  Then $B/J = C(E)$ for some closed 
subset $E$, and  now Holstynski's result is immediate 
using the well known correspondence between compact topological 
spaces and commutative $\cas$).  

4.  The known Banach-Stone theorem for surjective complete isometries between
$\cas$ is, formally, 
easy to derive from \ref{ma}.      
That is, if $T$ is also surjective then it is easy to see
that the *-homomorphism in \ref{ma}
maps onto $B$, and that $u$ is a unitary of $M(B)$.  However this
result also follows from much more elementary arguments.

\begin{corollary} \label{uca}
Let $T : A \rightarrow B$ be a unital completely contractive linear map
between unital $C^*$-algebras.  Then saying that 
 $T$ is a complete isometry is equivalent 
to any one of (ii)-(vi) in the previous theorem, but with the 
following changes: One may omit all mention of $u$, change 
`triple morphism' to `*-homomorphism',
and $K$ to $J^*$, in (iii) and (iv),  
and note that the expression $T(a) (1-p)$ in (ii) and (iv)
now coincides with $(1-p) T(a) (1-p)$.   
Indeed 
one may change the 
displayed equation in (vii) to read
$q_I(T(a)) = \pi(a)$ for all $a \in A$;
and $Z = D$ and $N = I$ in this case.     
\end{corollary}

We leave the proof of this as an exercise.  

 
We remark that Corollary \ref{uca} may be generalized to 
the case of completely isometric linear maps
$T : A \rightarrow B$ such that $T(1_A)$ is a unitary of  
the $\ca$ $B$.   For in this case, if $u = T(1_A)$, then 
$T' = u^* T(\cdot)$ is a unital complete isometry
$A \rightarrow B$.   One may then 
apply \ref{uca}  to obtain a nice form for  $T'$, and then 
use the fact that $T = u T'(\cdot)$.    For example, one 
obtains the existence of a projection $p \in B^{**}$ 
such that $T'(\cdot) (1-p)$ is a 1-1 *-homomorphism 
$\pi : A \rightarrow B^{**}$, then 
$\widehat{T(\cdot)} (1-p) = \hat{u} \pi(\cdot)$, with $\hat{u}$ 
unitary of $B^{**}$.    Clearly also we may write
$T(\cdot) (1-p) = v \pi(\cdot)$, where $v$ is a
unitary of the corner $(1-q) B^{**} (1-p)$ 
of $B^{**}$; namely $v = \hat{u} (1-p)$.    

More generally still we have:

\begin{proposition} \label{unin}  Suppose that $T :
A \rightarrow B$ is a complete isometry between unital
$\cas$, and suppose that  $T(A)$ contains a unitary (resp. isometry,
coisometry) of $B$.   Let $p, q$ be the associated projections
in $B^{**}$ as in \ref{ma} (see Remark 3 after 
\ref{ma}).  Then $U = \widehat{T(1)} (1-p)$
is a unitary (resp. isometry,
coisometry) of the corner $(1-q) B^{**} (1-p)$ of $B^{**}$,
and $\widehat{T(\cdot)} (1-p) = U \pi(\cdot)$ on $A$, where 
$\pi$ is the 1-1 *-homomorphism $A \rightarrow (1-p) B^{**} (1-p)$ 
considered in (v) of \ref{ma}.  
\end{proposition} 

\begin{proof}  As in \ref{ma} we have that 
$\widehat{T(\cdot)} (1-p)$ is a 1-1 triple morphism
$A \rightarrow B^{**} (1-p)$.    Let $Z$ be the triple 
system 
$\widehat{T(A)}(1-p) = 
(1-q) \widehat{T(A)} = (1-q) \widehat{T(A)} (1-p)$.   By \ref{maw},
$U$ is a unitary of $Z$.  If
$v = T(a_0)$ is an isometry (resp.  coisometry) of $B$, then $
V = \hat{v} (1-p)
= (1-q) \hat{v}$
is an isometry (resp.  coisometry)  of $(1-q) B^{**} (1-p)$ clearly.
Therefore $1-p$ is the identity of $Z^* Z$
if $v$ is an isometry, so that in this case $U^* U = (1-p)$.
Similarly in the `coisometry' case.
\end{proof}      

Remark.   If $T(A)$ contains a unitary $v$ of $B$, and if we let
$T' = v^* T(\cdot)$, then the triple system $Z'$ in $B$ 
generated by $T'(A)$ is a $\ca$ equal to $v^* Z$,
where $Z$ is the triple system generated by $T(A)$ as in
the proof above.  As in the 
proof of \ref{ma} we obtain 
a triple morphism $\theta' : Z' \rightarrow A$, and 
$N' = $ Ker $ \theta'$ is easily seen to equal
$v^* N$.  Thus the support projection $p' \in B^{**}$ 
for $N'$$^* N'$, equals $p$.   We also have that the 
associated support projection $q$ for $N N^*$ (see Lemma
\ref{moce}) equals $v p v^*$,
since $v p v^* z = v v^* z p = z p$ for all $z \in Z$. 
Thus $q$ is unitarily equivalent to $p$.    

\vspace{3 mm}

Let us amplify item (iv) of Theorem \ref{ma} in the case that $T(1_A) =
1_B$.   Recall that the projection $p$
discussed in the theorem lives in $B^{**}$.
 Thus if $H_u$ is the 
Hilbert space of the universal representation of $B$,
 then $p$ may be viewed as a 
projection on $H_u$.   Let $K = $ Ran $p$.
We may write $T = (1-p) T(\cdot) (1-p) + p T(\cdot) p$.
Setting $S = p T(\cdot) p$ (a completely positive unital map into
$B(K^\perp)$)
and $\theta = (1-p) T(\cdot) (1-p)$ (a unital *-homomorphism 
into $B(K)$),
we have  
$T = \pi \oplus S.$
That is,

\begin{corollary} \label{comp} Let $T : A \rightarrow B$ be a 
unital 
linear map
between unital $C^*$-algebras.  Let $H_u$ be the 
Hilbert space of the universal representation of $B$,     
Then $T$ is a 
complete isometry if and only if there exists 
a projection $p \in B''$ with respect to which $H_u$ splits as
$L \oplus L^\perp$, and with respect to which $T$ may be written 
as  $$T(a) = \left[ \begin{array}{ccl}
\pi(a) & 0 \\ 0 & S(a) \end{array} \right],$$
for a unital 1-1 *-homomorphism $\pi : A  \rightarrow (1-p) B''
(1-p) \subset B(L)$,
and for a completely positive unital 
$S : A  \rightarrow pB''p \subset B(L^\perp)$. 
Thus a compression of $T$  to the reducing subspace $L$
 of $H_u$,
is a 1-1 *-homomorphism. \end{corollary}

Of course a similar statement holds in the nonunital case.  

In the previous corollaries, and in \ref{ma},
  one may not hope to have $p \in B$ in general, even if 
$B$ is a $\wa$.    Indeed if $B \subset B(H)$, one cannot
hope to have the important projection $p$ above
in $B(H)$ (unless $H = H_u$ as we just discussed).   See 
Example \ref{nicex}.
Of course if $B$ is finite dimensional, then $p \in B$:    

\begin{lemma}  \label{trfrma}  Suppose that  $T : M_{m,n}
\rightarrow M_{r,s}$ is a triple morphism which 
is not identically zero.  Then 
$T$ is 1-1, and there exist unitary matrices $U$ and
$V$ of appropriat sizes, such that 
$T(A) = U diag \{ A , A , \cdots , A , 0 \} V$.
\end{lemma} 

\begin{proof}   This is essentially elementary, thus we 
merely sketch the proof.  
By \cite{Ham3}
2.1 (iv), we may view $T$ as the 1-2-corner of a
*-homomorphism $\pi : M_{m+n,m+n} \rightarrow 
M_{r+s,r+s}$.  Since $M_{m+n,m+n}$ is simple,
$\pi$ (and therefore also $T$) is  
1-1.   It is well known that such 
*-homomorphisms are unitarily equivalent to 
a map of the form $B \mapsto (B \otimes I_k)
 \; \oplus \; 0$.   Indeed looking at the 
elementary proof of this latter fact 
gives the desired result,
on inspection of the `1-2-corner'. 
\end{proof} 
 
\begin{corollary} \label{uc}    Suppose that $T : M_{m,n}
\rightarrow M_{r,s}$ is a linear map.   Then $T$ is a
complete isometry if and only if there exist unitaries 
$U$ and $V$ of appropriate sizes, such that $$UT(x)V =
\left[ \begin{array}{ccl}
x & 0 \\ 0 & S(x) \end{array} \right],$$
where $S : M_{m,n} \rightarrow M_{r-m,s-n}$ is a
complete contraction.   
\end{corollary}   

\begin{proof}   We apply \ref{bsttr}
to see that there exist projections
$e$ and $f$ in $M_s$ and $M_r$ respectively, such that
$$(1-f) T(z) = T(z) (1-e) = (1-f)
T(z) (1-e)$$
for all $z \in M_{m,n}$, and such that that these expressions
define a 1-1 triple morphism $\theta$ from $M_{m,n}$ into 
$M_{r,s}$.
That is, a reducing compression of $T$ is
a triple morphism. 
Now apply the last lemma to $\theta$.   
\end{proof}   

This corollary generalizes a result of Seth Hain 
(the case $n=1, m = r = s$) obtained  by different 
methods  during a Research Experiences for Undergraduates project 
in the summer of 2000, supervised by the first author
\cite{Hain}.  

The case that $m = n$ and $T(1) = 1$ is essentially 
7.1 in \cite{CEIOS}. 
   
It is easy to generalize the corollary to maps between
finite dimensional triple systems \cite{Hay}.

\begin{example} \label{nicex} \end{example}
In the light of the previous results, one might expect
that if $H$ is infinite dimensional,
and if $T : M_n \rightarrow B(H)$ is a complete isometry,
then there exists a projection $p$ in $B(H)$ with respect to
which $T(A) = U diag \{ A, S(A) \} V$ as in the previous
corollary.  The following example
shows that nothing like this is possible.

Suppose that $\{ \epsilon_k \}$ is a sequence of
numbers in $(0,1)$ strictly decreasing to $0$.
Define $\varphi^1_k \in (\ell^\infty_n)^*$ by
$\varphi^1_k((a_1,a_2, \cdots , a_n))  =
a_1 (1-\epsilon_k) + \frac{\epsilon_k}{n-1} \sum_{j=2}^n a_j$.
Note that $\varphi^1_k$ is a state of $\ell^\infty_n$.
Similarly define states $\varphi^r_k(\overset{\rightarrow}{a})
= a_r  (1-\epsilon_k) + \frac{\epsilon_k}{n-1} \sum_{j\neq r} a_j$.
Define  $\psi_k : \ell^\infty_n \rightarrow \ell^\infty_n$ by 
$\psi_k(\overset{\rightarrow}{a})
= (\varphi^1_k(\overset{\rightarrow}{a}), \cdots ,
\varphi^n_k(\overset{\rightarrow}{a}))$.
Then $\psi_k$ is a unital contraction, and hence is
completely contractive and completely positive.   Note also that
$$|\varphi^r_k(\overset{\rightarrow}{a})|
\geq |a_r| (1-\epsilon_k) - |\frac{\epsilon_k}{n-1} \sum_{j\neq r} a_j|
\geq |a_r| - 2 \epsilon_k \max_k |a_k|,$$
from which it is easy to see that
$$\Vert \psi_k(\overset{\rightarrow}{a}) \Vert_\infty
\; \geq \; (1 - 2 \epsilon_k) \max_k |a_k|
\; = \; (1 - 2 \epsilon_k) \Vert \overset{\rightarrow}{a} \Vert_\infty.$$
Thus $\sup_k \Vert \psi_k(\overset{\rightarrow}{a}) \Vert_\infty =
\Vert \overset{\rightarrow}{a} \Vert_\infty$, for all
$\overset{\rightarrow}{a} \in \ell^\infty_n$.

Define a map $\Phi : \ell^\infty_n \rightarrow \ell^\infty(\ell^\infty_n)$
by $\Phi(\overset{\rightarrow}{a}) = ( \psi_k(\overset{\rightarrow}{a}) )_k$.  
Then $\Phi$ is a unital isometry between subspaces of 
$C(K)$ spaces,
and is therefore a complete isometry (\cite{Pau} Theorem 3.8).
We consider the canonical `main diagonal'
*-representation $\pi : \ell^\infty(\ell^\infty_n) \rightarrow B(H)$, where
$H = \oplus_{k=1}^\infty H_i$, and 
each $H_i$ is a copy of $\mathbb{C}^n$.
Let $\Psi = \pi \circ \Phi$.  Then $\Psi$ is a unital complete isometry.
We claim that if there exists a projection $p \in B(H)$ such that
$$(1-p) \Psi(\cdot) = \Psi(\cdot) (1-p)$$ is a *-homomorphism,
then $p = 1$ and this *-homomorphism is the zero map.   
Here $\Psi : \ell^\infty_n \rightarrow B(H)$, but it is easy to see how
$\Psi$ may be extended to a unital complete isometry $M_n \rightarrow B(H)$
(extend the $\psi_k$ above to 
maps $M_n \rightarrow M_n$ by viewing 
$\ell^\infty_n$ as the `main diagonal' and multiplying off-diagonal 
terms by $1-\epsilon_k$).  Let $A = M_n$ or $\ell^\infty_n$.    

To this end, note that if there were such $p$, then
$p$ commutes with $D_i \overset{def}{=} 
\Psi(\overset{\rightarrow}{e_i})$ for each
$i = 1 , \cdots , n$.   Suppose that $A_{ij} = P_i p P_j$,
where $P_i$ is the projection of $H$ onto $H_i$.
Let $E_{ij}$ be the matrix units in $M_n$.
Note that
$\psi_k(\overset{\rightarrow}{e_1}) = ((1 - \epsilon_k),\frac{\epsilon_k}{n-1},
\cdots ,\frac{\epsilon_k}{n-1})$,
 so that $\Psi(\overset{\rightarrow}{e_1})$
is the matrix $$D_1 \; = \;
diag \{ (1 - \epsilon_1) , \frac{\epsilon_1}{n-1} I_{n-1} ,
(1 - \epsilon_2) , \frac{\epsilon_2}{n-1} I_{n-1} , \cdots \}.$$
That $p$ commutes with $D_1$ implies that
$A_{ij} = 0$ if $i \neq j$; and that
$A_{ii}$ commutes with $E_{11}$ for all $i$.   Since
$p$ is a projection we obtain
 $A_{ii} E_{11} = E_{11}$ or $A_{ii} E_{11} = 0$.
A similar argument applied to $D_r$ for $2 \leq r \leq n$ implies
that $A_{ii} E_{rr}$ equals $E_{rr}$ or $0$.  Thus
$p$ and $1-p$ are `diagonal matrices', with either $0$'s or $1$'s on the
diagonal.  Since  $(1-p) \Psi(\cdot)$
is a *-homomorphism, we know that
$(1-p) D_r$ is a projection for each $r$.   However then it is clear
from the form of $1-p$ and $D_r$, that $(1-p) D_r$
must be the zero projection for every $r$.   Also, $p = 1$.
Therefore $(1-p) \Psi(\cdot)$ is the zero map.

Indeed suppose that there exists projections $p, q \in B(H)$ such that
$$(1-q) \Psi(\cdot) = \Psi(\cdot) (1-p)$$ is a triple morphism
$\psi$.  Inserting 
$1$, gives $p = q$, and so $\psi(A)$ is a
triple system containing $(1-p)$.  Hence $\psi(A)$ is a unital $\ca$,
and $\psi$ is a unital *-homomorphism.  
This puts us back in the previous
case, so that $p = q = 1$ and $\psi$ is the zero map.

Acknowledgments:   We thank Louis Labuschagne for 
helpful discussions and input during a visit of the first author.
We are indebted to the University of South 
Africa for support for that visit.


\begin{thebibliography}{99}
\bibitem{Arv1}  W. B. Arveson,
{\em Subalgebras of }$C^{*}-${\em algebras,}
 Acta Math. {\bf 123 }(1969), 141-224.

 
\bibitem{Arv2}   W. B. Arveson,
{\em Subalgebras of }$C^{*}-${\em algebras II,}  Acta Math.
 {\bf 128} (1972), 271-308.
 
\bibitem{BShi}  D. P. Blecher,  {\em  The Shilov boundary
of an operator space, and the characterization
theorems,}  J. Funct. An. {\bf 182} (2001), 280-343.

\bibitem{BEZ}  D. P. Blecher, E. G. Effros and V. Zarikian,
{\em One-sided $M$-Ideals and Multipliers in Operator Spaces
I}, to appear Pacif. J. Math.

\bibitem{BL}  D. P. Blecher and L. E. Labuschagne,  
Work in progress.

\bibitem{CE}  M. D. Choi and E.G. Effros, {\em The completely 
positive lifting problem for $\cas$}, Ann. Math. {\bf 104}
(1976), 585-609.

\bibitem{CEIOS}  M. D. Choi 
and E.G. Effros, {\em  Injectivity and operator spaces,} J. Funct. An.
{\bf 24} (1977), 156-209.

\bibitem{Dixon}  P. G. Dixon, {\em Non-closed sums in 
closed ideals in Banach algebras,}
Proc. Amer. Math. Soc., {\bf 128} (2000), 3647-3654.

\bibitem{EOR}  E. G.
Effros, N. Ozawa and Z-J. Ruan, {\em On injectivity and
nuclearity for operator spaces,} Duke J. Math. {\bf 110} (2001),
489-521.   

\bibitem{Hain}  S. Hain, {\em Algebra and Matrix Normed Spaces,}
Rose-Hulman Undergrad. Math. J. {\bf 2} (2) (2001). 
 
\bibitem{Ham}  M. Hamana, {\em Injective
envelopes of operator systems,}
Publ. R.I.M.S. Kyoto Univ. {\bf 15} (1979), 773-785.
 
\bibitem{Ham3}  M. Hamana, {\em Triple  envelopes and Silov
boundaries of operator spaces,}  Math. J. Toyama University
{\bf 22} (1999), 77-93.

\bibitem{Hay}  D. Hay,  Notes for Ph. D. Thesis project supervised
by D. P. Blecher.   

 
\bibitem{Hol}  W.  Holsztynski, {\em Continuous mappings 
induced by isometries of spaces of continuous functions,
} Studia Math. {\bf 24} (1966), 133-136.
 
 
\bibitem{Ka}  R. V. Kadison, {\em Isometries of
operator algebras,}  Ann. of Math. {\bf 54} (1951),
325-338.

\bibitem{Ki}  E. Kirchberg, {\em On restricted peturbations
in inverse images and a description of normalizer algebras
in $C^*$-algebras,}  J. Funct. An. {\bf 129} (1995), 1-34.

\bibitem{Ki2}  E. Kirchberg, {\em Commutants of unitaries
in UHF algebras and functorial properties of 
exactness,}  J. Reine Angew Math. {\bf 452} (1994), 39-77.

 \bibitem{KW}  E. Kirchberg and S. Wassermann, {\em
C$^*-$algebras generated by operator systems,}
J. Funct. Analysis {\bf 155} (1998), 324-351.

 
\bibitem{Mg}  B. Magajna, {\em Hilbert modules and tensor products
of operator spaces,} Linear Operators, Banach Center Publ. Vol. 38,
Inst. of Math. Polish Acad. Sci. (1997), 227-246.

\bibitem{Pau} V. I. Paulsen, {\em
Completely bounded maps and dilations,} Pitman Research
Notes in Math., Longman, London, 1986.

\bibitem{Ped} G. Pedersen, {\em $C^*$-algebras
and their automorphism groups,}
Academic Press (1979).

\bibitem{Ri2}  M. A. Rieffel, {\em Morita equivalence for operator
algebras, }Proceedings of Symposia in Pure Mathematics
 {\bf 38} Part 1
(1982), 285-298.

\bibitem{Rif2}  M. A. Rieffel, {\em  Unitary representations
of group extensions; an algebraic approach to the theory
of Mackey and Blattner,}  Studies in Analysis, Adv. Math.
Suppl. Stud. {\bf 4} (1979), 43-82.
 
\bibitem{Rieffel2}  M. A. Rieffel,
{\em Morita equivalence for }$C^{*}-${\em algebras and
W}$^{*}-${\em algebras,} J. Pure Appl. Algebra {\bf 5 }(1974),
51-96.
 
\bibitem{Wi}  G. Wittstock, {\em Extensions of
completely bounded module
morphisms}, Proceedings of conference on operator algebras and
group representations, Neptum, 238-250,
Pitman (1983).
\end{thebibliography}
\end{document}